\documentclass[11pt,reqno]{amsart}
\usepackage{amssymb,amsmath,amsthm,amsxtra,calc,bm, color}
\usepackage{enumerate}
\usepackage[margin=.8in]{geometry}
   \usepackage{stmaryrd}
\usepackage[scr=boondoxo]{mathalfa}
\usepackage{mathtools, makecell}
\usepackage{tabularx}
\usepackage[labelsep=period]{caption}
\usepackage{float, color, multirow}
\usepackage[labelsep=period]{caption}
\usepackage{todonotes, subcaption}
\usepackage{float}
\DeclareMathOperator{\Gal}{Gal}
\DeclareMathOperator{\GL}{GL}
\DeclareMathOperator{\SL}{SL}
\DeclareMathOperator{\PGL}{PGL}
\usepackage{makecell}

\usepackage{relsize}
\usepackage{mleftright}
\mleftright

\usepackage{enumitem}

\usepackage{nccmath}
\usepackage{cases}
\usepackage{graphicx, standalone}
\usepackage{hyperref, commath}
\usepackage[final]{microtype}
\usetikzlibrary{arrows.meta}

\newcommand{\ord}{\operatorname{ord}}

\def\Z{\mathbb{Z}}
\def\Q{\mathbb{Q}}

\def\H{\mathbb{H}}

\def\C{\mathbb{C}}
\def\P{\mathbb{P}}
\def\F{\mathbb{F}}

\def\SL{{\rm SL}}

\def\GL{{\rm GL}}

\usepackage{xspace}

\renewcommand{\tilde}{\widetilde}

\DeclareMathOperator{\new}{new}

\newtheorem{theorem}{Theorem}[section]
\newtheorem{lemma}[theorem]{Lemma}

\newtheorem{proposition}[theorem]{Proposition}
\newtheorem{conjecture}[theorem]{Conjecture}

\theoremstyle{remark}

\numberwithin{equation}{section}

\mathtoolsset{showonlyrefs}

\makeatletter
\newcommand*\bigcdot{\mathpalette\bigcdot@{.7}}
\newcommand*\bigcdot@[2]{\mathbin{\vcenter{\hbox{\scalebox{#2}{$\m@th#1\bullet$}}}}}
\makeatother

\usepackage{listings}

\lstdefinestyle{myListingStyle} 
    {
        basicstyle = \small\ttfamily,
        breaklines = true,
    }
\usepackage[nohyperlinks]{acronym}
\renewcommand{\url}[1]{#1}
\makeatletter
\newcommand*{\rom}[1]{\expandafter\@slowromancap\romannumeral #1@}
\makeatother

\makeatletter
\def\imod#1{\allowbreak\mkern5mu({\operator@font mod}\,\,#1)}
\makeatother
\allowdisplaybreaks


\makeatletter
\@namedef{subjclassname@2020}{%
  \textup{2020} Mathematics Subject Classification}
\makeatother

\begin{document}
\title[]{Exceptional Congruences for Eta-quotient newforms}
%\author{Eddie O'Sullivan, Henry Stone, Swati, and Xiaolan Jin}

\begin{abstract}
 In 1973, Swinnerton-Dyer completely classified all congruences for coefficients of normalized eigenforms in weights $k \in \{12, 16, 18, 20, 22, 26\}$ on $\Gamma_{0}(1) = \textup{SL}_{2}(\mathbb{Z})$ using the theory of modular Galois representations. In this paper, we classify congruences of Type I and Type II considered by Swinnerton-Dyer for the coefficients of eta-quotient newforms in $S_{k}(N, \chi)$. When $k \geq 2$, we prove them using the theory of modular forms modulo primes. We also prove extensions of these congruences modulo prime powers.
\end{abstract}

\author{Eddie O'Sullivan}
\address{Department of Mathematics, Colby College, Waterville, ME 04901, USA}
\email{emosul26@colby.edu}

\author{Henry Stone}
\address{Department of Mathematics, University of Michigan, Ann Arbor, ME 04901, USA}
\email{henstone@umich.edu}

\author{Swati}
\address{Department of Mathematics, University of South Carolina, Columbia, SC 29208, USA}
\email{s10@email.sc.edu}

\author{Xiaolan Jin}
\address{Department of Mathematics, University of Illinois Urbana-Champaign \\ Urbana, IL 61801, USA}
\email{seanj4@illinois.edu}

\date{}


\subjclass[2020]{}

\maketitle

\section{Introduction and Statement of Results}
\subsection{Introduction}
For $z \in \H$, we consider the Ramanujan Delta function $\Delta(z)$ given by the following infinite product:
 \[
 \Delta(z) = \sum_{n \geq 1} \tau(n) q^n = q \prod_{n \geq 1} (1 - q^n)^{24}, \: \text{where} \: q = e^{2 \pi i z}.
 \]
 For primes $\ell \in \{2,3,5,7,23, 691\}$ and $m \geq 1$, Ramanujan \cite{ramanujan1916certain}, Wilton \cite{wilton1930congruence}, Kolberg \cite{kolberg1962congruences}, Ashworth \cite{ashworth1968congruence}, and Lehmer \cite{lehmer8note}  proved congruences of the form 
 \begin{equation} \label{cong_delta}
 \tau(n) \equiv c \cdot \sigma_{\nu}(n) \pmod{\ell^m} \: \text{for some} \: c \in \Z.
 \end{equation}
Here, $n$ is a positive integer or $n \equiv r \pmod{\ell^{m'}}$ where $m' \mid m$ and $\sigma_{\nu}(n) = \sum_{d \mid n} d^{\nu}$ for a non-negative integer $\nu$. In 1973, Swinnerton-Dyer \cite{swinnerton1973ℓ} provided a general framework explaining every congruence for the coefficients of a normalized eigenform on $\SL_2(\Z)$ with integer coefficients. He used the theory of modular Galois representations developed by Deligne \cite{Deligne1968-1969} for $k \geq 2$.

\begin{theorem}[Serre-Deligne]
   Let $k \geq 2$ be an integer. Suppose that
$f = \sum_{n \geq 1} a_{f}(n) q^{n} \in S_{k}(\Gamma_{0}(N), \chi)$ 
is a normalized eigenform with $\ell$-adic integer coefficients. Then there exists a continuous homomorphism
\begin{equation}
    \rho_{\ell, f} : \Gal\left( \overline{\Q}/\Q\right) \longrightarrow \GL_{2}(\Z_{\ell})
\end{equation}
for each prime $\ell$, depending on $f$, such that $\rho_{\ell, f}(Frob_{p})$ satisfies the characteristic polynomial
\[
X^2 - a_{f}(p)X + \chi(p) p^{k - 1},
\]
for each $p \neq N \ell$. 
\end{theorem}

When $k=2$ and $f \in S_{2}(\Gamma_{0}(N))$ has integer Fourier coefficients, then Eichler-Shimura implies that there exists an elliptic curve $E/\Q$ of conductor $N$ such that $L(f,s) = L(E,s)$. Further, for every prime $p \nmid N$, the trace of Frobenius of $E$ satisfies
\begin{equation}\label{trace}
a_E(p) = p+1 - \#E(\F_p),
\end{equation}
and one has the identity
\begin{equation} \label{Eichler-Shimura}
a(p) = a_E(p).
\end{equation}
 
\medskip
\noindent
We denote by 
\[
\tilde{\rho}_{\ell, f} : \Gal\left(\overline{\Q}/\Q\right) \longrightarrow \GL_{2}(\F_{\ell}).
\]
the reduction of $\rho_{\ell, f}$ modulo $\ell$.
We note that Swinnerton-Dyer calls the primes $\ell$ appearing in congruence \ref{cong_delta} \textit{exceptional}. Let $\ell > 3$ be a prime. We say $\ell$ is exceptional for $f$ if and only if the image of  $\tilde{\rho}_{\ell, f}$ doesn't contain $\SL_2(\F_{\ell})$. When $\ell \in \{2,3\}$, the condition is sufficient for $\ell$ to be exceptional for $f$.

Let $G$ denote the image of $\tilde{\rho}_{\ell, f}$ in $\GL_{2}(\F_{\ell})$ and let $H$ denote the image of $G$ in $\PGL_{2}(\F_{\ell})$ under the natural projection. The exceptional primes $\ell$ are classified into the following types, and in each case the Fourier coefficients of $f$ satisfy specific congruences, as described in the following lemma.

\begin{lemma}
Let $f$, $G$, and $H$ be as defined above. Suppose that $\ell$ is exceptional for $f$. Then exactly one of the following occurs:
\begin{enumerate}
    \item If $G$ is contained in a Borel subgroup of $\GL_{2}(\F_{\ell})$, then $\ell$ is an exceptional prime of type I for $f$. In this case, there exists an integer $m$ such that 
    \begin{equation*}
        a(n) \equiv n^{m}  \sigma_{k - 1 - 2m}(n) \pmod{\ell} 
        \: \text{for all} \;\; n \geq 1  \: \text{with} \: \gcd(n, \ell) = 1.
    \end{equation*}
    \item If $H$ is dihedral, then $\ell$ is an exceptional prime of type II for $f$. In this case, we have
      \begin{equation*}
        a(n) \equiv 0 \pmod{\ell} 
        \: \text{for all } \: n \geq 1 \: \text{ with } \: \left( \tfrac{n}{\ell} \right) = -1. 
      \end{equation*}
    \item If $H \cong A_4, S_{4}$ or $A_5$, then $\ell$ is an exceptional prime of type III for $f$. When $H \cong S_4$, for every prime $p \nmid N \ell$, we have
    \begin{equation} \label{image_S4}
        \frac{a(p)^{2}}{\chi(p) p^{k - 1}} \equiv 0, 1, 2, \; \text{or} \:\: 4 \pmod{\ell}.
    \end{equation}
When $H \cong A_4$, the set \{0,1,2,4\} in \eqref{image_S4} is replaced by \{0,1,4\}. 
When $H \cong A_5$, either the set \{0,1,2,4\} in \eqref{image_S4} is replaced by \{0,1,4\} or we have
\[
a(p)^4 - 3 \chi(p) p^{k - 1} a(p)^2 + \chi(p)^2 p^{2k - 2} \equiv 0 \pmod{\ell}
\]
          \end{enumerate}
\end{lemma} 
\noindent
We begin by recalling several important results concerning the classification of exceptional primes and congruences between modular forms. Fundamental contributions in this area are due to Swinnerton-Dyer \cite{swinnerton1973ℓ, swinnerton2006ℓ}, Serre \cite{serre2006congruences}, Ribet \cite{Ribet1975, ribet1985adic}, and Momose \cite{momose1981adic}, who showed that only finitely many primes $\ell$ are exceptional for a given eigenform $f$. The foundational paper \cite{swinnerton1973ℓ} classified exceptional congruences of Type I and Type II for coefficients of normalized eigenforms with integer coefficients on $\SL_{2}(\Z)$, and conjectured a unique congruence of Type III for this family.
\begin{conjecture}[Serre-Swinnerton-Dyer] \label{swd_59}
The prime $\ell = 59$ is an exceptional prime of type III for 
\[
f(z) = E_{4}(z)\,\Delta(z) \in S_{16}(\Gamma_{0}(1)).
\]
\end{conjecture}

In 1983, Haberland \cite{haberland1983perioden} proved the Conjecture \ref{swd_59} using techniques from Galois cohomology. 
\medskip

\noindent
We shift our attention to eta-quotients. The Dedekind eta-function is a function of the complex upper half-plane $\mathbb{H}$ defined by
\begin{equation}
\label{eta1}
\eta(z) = q^{\frac{1}{24}}\prod_{n = 1}^{\infty}(1 - q^n), \ \ q = e^{2\pi iz}.
\end{equation}
The eta-function is a building block for modular forms. An eta-quotient of level $N$ is a function of the form
\begin{equation}
\label{etaquot}
f(z) = \prod_{\delta \mid N} \eta(\delta z)^{r_{\delta}},
\end{equation}
where $\delta$ and $r_{\delta}$ are integers with $\delta \geq 1$. 
The aim of this paper is to complete a similar classification for coefficients of eta-quotient newforms. In 1996, Y. Martin \cite{martin1996multiplicative} proved that only finitely many such eta-quotients exist and provided a complete list (see Table \rom{1} in Section 7). In 2003, Boylan \cite{boylan2003exceptional} proved all congruences of type III for the coefficients of eta-quotient newforms. More generally, he provided a framework to obtain every congruence for $f$ arising from odd, complex, two-dimensional Galois representations with projective image isomorphic to $S_{4}$. In addition, he gave an alternate proof of the congruence modulo $59$ as stated in \eqref{image_S4} with $k = 16$ and primes $p \neq 59$.

\begin{theorem}[Boylan]\label{boylan_type3}
Consider the following eta-quotient newforms 
\begin{align*}
f_1(z) &= \eta^{4}(2z)\,\eta^{4}(4z) \in S_{4}(\Gamma_{0}(8)),  
& f_2(z) &= \frac{\eta^{16}(4z)}{\eta^4(2z) \eta^4(8z)} = f_1(z) \otimes \left( \tfrac{-4}{\cdot} \right) \in S_{4}(\Gamma_{0}(128)), \\ 
f_3(z) &= \eta^{12}(2z) \in S_{6}(\Gamma_{0}(4)),  
& f_4(z) &= \frac{\eta^{36}(4z)}{\eta^{12}(2z) \eta^{12}(8z)} = f_3(z) \otimes \left( \tfrac{-4}{\cdot} \right) \in S_{6}(\Gamma_{0}(64)).
\end{align*}
Then we have
\begin{enumerate}
    \item $\ell = 11$ is an exceptional prime of type III for $f_1(z)$ and $f_2(z)$. 
    \item $\ell = 19$ is an exceptional prime of type III for $f_3(z)$ and $f_4(z)$. 
\end{enumerate}
\end{theorem}

In 2000, Kiming and Verrill \cite{kiming2005modular} employed a different approach to prove the congruences stated in Theorem~\ref{boylan_type3}, as well as the congruence modulo $59$.

\subsection{Statement of Results}
We now turn to our work. Let $f = \displaystyle{\sum_{n \geq 1}} a_{f}(n) q^{n} \in S_{k}(\Gamma_{0}(N), \chi) \cap \Z_{(\ell)}[[q]]$ be an eta-quotient newform. We require some terminology which generalizes the terminology in \cite{swinnerton1973ℓ}. We say that $f$ satisfies a Type I congruence precisely when there exists real Dirichlet characters $\psi, \phi$ modulo $N$ and non-negative integers $m' \geq m$ such that for all primes $p \nmid N \ell$, we have
     \begin{equation}
        \label{Trcond}
        a(p) \equiv \psi(p) p^{m} + \phi(p) p^{m'} \pmod \ell
    \end{equation}
    and
    \begin{equation}
        \label{Detcond}
        \psi(p)\phi(p)p^{m' + m} \equiv \chi(p) p^{k-1} \pmod \ell.
    \end{equation}

Our first theorem provides a classification of primes $\ell$ that are exceptional of Type I for an eta-quotient newform $f$. We note that for a Type I congruence, we have $\psi = \phi$ (see \S 3.1 for details).
\begin{theorem} \label{type1_cong}
Let $f$ and $\psi$ be as defined above. For even $k \geq 4$, let $G_{k}(z) = \frac{-B_k}{2k} + \sum_{n \geq 1} \sigma_{k - 1}(n) q^n$ and $E_k(z) = \frac{-2k}{B_k} G_{k}$
denote Eisenstein series of weight $k$. For $k = 2$ and $N \geq 2$, we replace $G_k$ by $E_{2, N}$ as defined in \eqref{defn_E2N}.
Let $\psi$ be a real-valued Dirichlet character modulo $N$, and let $1_{N}$ denote the trivial character modulo $N$.
\begin{enumerate}
    \item If $3 \leq m' - m + 1 \leq \ell - 2$, then we have
    \begin{equation}
        \theta(f \otimes 1_N) \equiv \theta^{m + 1} (G_{m' - m + 1} \otimes \psi1_N) \pmod{\ell}.
    \end{equation}
    Moreover, we have $\ell < k$ or $\ell \mid (a(p) - \psi(p) \: \sigma_{k - 1}(p))$ for all primes $p \nmid N$.
    \item If $m' - m + 1 \in \{ 2, \ell - 1 \} $, then we have
    \begin{equation}
        \theta(f \otimes 1_N) \equiv \theta^{m + 1} (G_{\ell + 1} \otimes \psi1_N) \pmod{\ell}.
    \end{equation}
Further, we have $\ell < k$.
\end{enumerate}

  \begin{table}[H]
  \centering
  \caption{Values of parameters in Theorem \ref{type1_cong}}
  \begin{subfigure}{0.48\textwidth}
    \centering
    \begin{tabular}{|c|c|c|c|c|c|}
        \hline
        $f(z)$ & \multicolumn{5}{c|}{type \rom{1}} \\
        \cline{2-6}
        & $\ell$ & $m$ & $m'$  & $\psi$ & $(k,N)$ \\
        \hline
        \makecell{$\Delta(z)$} 
        & \makecell{3 \\ 5 \\ 7 \\ 691} 
        & \makecell{0 \\ 1 \\ 1 \\ 0} 
        & \makecell{1 \\ 2 \\ 4 \\ 11}
        & $1_{1}$ & $(12,1)$ \\
        \hline
        \makecell{$\eta(z)^8 \eta(2z)^8$} 
        & \makecell{2 \\ 3 \\ 5 \\ 17} 
        & \makecell{0 \\ 0 \\ 1 \\ 0} 
        & \makecell{1 \\ 1 \\ 2 \\ 7}
        & $1_{2}$ & $(8, 2)$ \\
        \hline
        \makecell{$\eta(z)^6 \eta(3z)^6$} 
        & \makecell{2 \\ 3 \\ 13} 
        & \makecell{0 \\ 0 \\ 0} 
        & \makecell{1 \\ 1 \\ 5}
        & \makecell{ $1_{3}$ \\ $\left( \frac{\cdot}{3}\right)$\\ $1_{3}$} & $(6, 3)$ \\
        \hline
        \makecell{$\eta(2z)^{12}$}
        & \makecell{2 \\ 3}
        & \makecell{0 \\ 0} 
        & \makecell{1 \\ 1}
        & \makecell{$1_{4}$ \\ $1_{4}$} & $(6, 4)$ \\
        \hline
        \makecell{$\eta(z)^4 \eta(5z)^4$} 
        & \makecell{2 \\ 5 \\ 13} 
        & \makecell{0 \\ 0 \\ 0} 
        & \makecell{1 \\ 3\\ 3}
        & \makecell{$1_{5}$ \\ $\left( \frac{\cdot}{5}\right)$\\ $1_{5}$} & $(4, 5)$ \\
        \hline
    \end{tabular}
    %\subcaption{Part 1}
    \label{tab:type1_cong}
  \end{subfigure}
  \hfill
  \begin{subfigure}{0.48\textwidth}
    \centering
    \small
    \begin{tabular}{|c|c|c|c|c|c|}
        \hline
       $f(z)$ & $\ell$ & $m$ & $m'$  & $\psi$ & $(k, N)$ \\
        \hline
        \makecell{$\eta(z)^2 \eta(2z)^2 \eta(3z)^2 \eta(6z)^2$} 
        & \makecell{2 \\ 3 \\ 5} 
        & \makecell{0  \\ 0 \\ 0} 
        & \makecell{1 \\ 1 \\ 3}
        & \makecell{$1_{6}$ \\ $1_{2} \left( \frac{\cdot}{3}\right)$\\ $1_{6}$} & $(4, 6)$\\
        \hline
        \makecell{$\eta(3z)^8$} 
        & \makecell{2 \\ 3} 
        & \makecell{0 \\ 0} 
        & \makecell{1 \\ 1}
        & \makecell{$1_{9}$ \\ $1_{3} \left( \frac{\cdot}{3}\right)$} & $(4, 9)$\\ 
        \hline
        \makecell{$\eta(z)^2 \eta(11z)^2$} 
        & \makecell{5} 
        & \makecell{0} 
        & \makecell{1}
        & $1_{11}$ & $(2, 11)$\\ 
        \hline
        \makecell{$\eta(z) \eta(2z) \eta(7z) \eta(14z)$} 
        & \makecell{2 \\ 3} 
        & \makecell{0 \\ 0} 
        & \makecell{1 \\ 1}
        & $1_{14}$ & $(2, 14)$\\ 
        \hline
         \makecell{$\eta(z) \eta(3z) \eta(5z) \eta(15z)$} 
        & \makecell{2} 
        & \makecell{0} 
        & \makecell{1}
        & $1_{15}$ & $(2, 15)$\\ 
        \hline
        \makecell{$\eta(2z)^2 \eta(10z)^2$} 
        & \makecell{2 \\ 3} 
        & \makecell{0 \\ 0} 
        & \makecell{1 \\ 1}
        & $1_{20}$ & $(2, 20)$ \\ 
        \hline
        \makecell{$\eta(3z)^2 \eta(9z)^2$} 
        & \makecell{3} 
        & \makecell{0} 
        & \makecell{1}
        & \makecell{$1_{9} \left( \frac{\cdot}{3} \right)$} & $(2, 27)$\\ 
        \hline
        \makecell{$\eta(6z)^4$} 
        & \makecell{2 \\ 3} 
        & \makecell{0  \\ 0} 
        & \makecell{1 \\ 1}
        & \makecell{$1_{36}$ \\ $1_{12} \left( \frac{\cdot}{3} \right)$} & $(2, 36)$\\ 
        \hline
    \end{tabular}
   % \subcaption{Part 2}
    \label{tab:continuation}
  \end{subfigure}
\end{table}
\end{theorem}
%\newpage
\textbf{Remarks:} 
    \begin{enumerate}
    \item There are no congruences of type I for eta-quotient newforms with non-trivial character $\chi$. 
    \item The eta-quotients $ \eta(3z)^8$ and $\eta(6z)^4$ are CM newforms.
    \item There are exceptional primes of Type I in the theorem that divide the level $N$ of the eta-quotient newform $f(z)$. Using $\eta(\ell z) \equiv \eta(z)^{\ell} \pmod{\ell}$ which follows from the Binomial Theorem mod $\ell$, we see that these primes are exceptional of Type I for an eta-quotient newform $g(z) \equiv f(z) \pmod{\ell}$ with $\ell$-free level. 
        \item When $k = 2$, congruences for Hecke eigenvalues translate directly into congruences for point counts on the associated elliptic curve. For instance, we can associate the elliptic curve $E: 
    y^2 - y = x^3 - x$, to the eta-quotient newform 
   $ f(z) = \eta(z)^2 \eta(11z)^2$,
    whose coefficients obey the congruence
    $a(p) \equiv p + 1 \pmod 5$,
    for $p \nmid 55$. Using \eqref{trace} and \eqref{Eichler-Shimura}, we deduce that
    $p+1 - \#E(\mathbb{F}_p) \equiv p +1 \pmod 5$
   implying that $\#E(\mathbb{F}_p) \equiv 0 \pmod 5$.
        \item The above congruences can be interpreted combinatorially. For instance, let 
        \[
        \sum_{n \geq 0}v(n)q^n = \prod_{n \geq 1} (1 - q^{ n})^{2} (1 - q^{ 11 n})^{2}.
        \]
        Then we have $v(n) = v_{e}(n) - v_{o}(n)$, where $v_{e}(n)$ and $v_{o}(n)$ denote the number of 4-colored partitions of $n$ into an even (resp. odd) number of dictinct parts where the parts of last 2 colors are multiples of 11. Hence, we obtain for all $(n, 11) = 1$, 
        \[
        v(n - 1) \equiv a_{11}(n) \pmod{5} \:\: \text{where} \:\: a_{11}(n) = \sum_{\substack{d \mid n \\ \gcd(d, 11) = 1}} d.
        \]
    \end{enumerate}
    \par
Let $\ell \geq $3 be prime. We say that $f$ satisfies a Type II congruence modulo $\ell$ precisely when for all primes $p \nmid N \ell$, we have 
\begin{equation}\label{type2}
    a(p) \equiv 0 \pmod{\ell} \: \text{whenever} \: \left( \frac{p}{\ell} \right) = -1.
\end{equation}
We now classify exceptional congruences of Type II for eta-quotient newforms.

\begin{theorem} \label{type2_cong}
      Let $f \in S_{k}(\Gamma_0(N), \chi)$ be an eta-quotient newform and let $\ell$ be an exceptional prime of Type II for $f$. Then we have
      \[
      \Theta^{\frac{\ell + 1}{2}}(f \otimes 1_N) \equiv \Theta(f \otimes 1_N) \pmod{\ell},  
      \]
      and we have $\ell < k$ or
     \[
    \ell = \begin{cases}
        2k - 1 \quad \text{if} \quad f \mid U_{\ell} \not\equiv 0 \pmod{\ell} \\
        2k - 3 \quad \text{if} \quad f \mid U_{\ell} \equiv 0 \pmod{\ell}
    \end{cases} \; \text{where} \; U_{\ell} \; \text{is as defined in \eqref{UVq}}.
    \]
      \end{theorem}

      \vspace{-6mm}
\renewcommand{\arraystretch}{1.5} % Adjust the vertical spacing
\setlength{\tabcolsep}{12pt} % Adjust the horizontal spacing
\begin{table}[ht]
  \centering
  \small
  %\resizebox{1.3\textwidth}{!}{%
 % \resizebox{\textwidth}{!}{%
\begin{tabular}{|c|c|c|c|}
\hline
$f(z)$ & $\ell$ & $(k, N)$ & $\chi$\\ 
\hline
$\Delta(z)$ & 23 & $(12, 1)$ & $1_{1}$ \\
 \hline
 $\eta(z)^8 \eta(2z)^8$ & 3  & $(8, 2)$ & $1_{2}$\\
 \hline
$\eta(z)^6 \eta(3z)^6$ & 3  & $(6, 3)$ & $1_{3}$\\
 \hline
\makecell{$\begin{aligned}\eta(2z)^{12}, \: &\frac{\eta(4z)^{36}}{\eta(2z)^{12} \eta(8z)^{12}} \\ &= \eta(2z)^{12} \otimes \left(\frac{-4}{\cdot}\right)\end{aligned}$} & 3, 11  & $(6, 4), (6, 16)$ & $1_2, 1_2$ \\
 \hline
 \makecell{$\begin{aligned}\eta(z)^{4} \eta(2z)^2 \eta(4z)^4,  \:&\frac{\eta(8z)^{38}}{\eta(4z)^{14} \eta(16z)^{14}}\\ &= \eta(z)^{4} \eta(2z)^2 \eta(4z)^4 \otimes \left( \frac{-8}{\cdot} \right) \end{aligned}$} & 7 & (5, 4), (5, 64) & $\left( \frac{-4}{\cdot} \right), \left( \frac{-4}{\cdot} \right)$\\
 \hline
$\eta(z)^{3} \eta(7z)^3$ & 3  & $(3, 7)$ & $\left( \frac{-7}{\cdot} \right)$ \\
 \hline
$\eta(3z)^{8}$ &3,  5, 7 &  $(4, 9)$ & $1_3$\\
 \hline
\makecell{$\begin{aligned}\eta(2z)^{3} \eta(6z)^3,  \:&\frac{\eta(4z)^{9} \eta(12z)^9}{\eta(2z)^{3} \eta(6z)^3 \eta(8z)^3 \eta(24z)^3} \\ &= \eta(2z)^{3} \eta(6z)^3 \otimes \left(\frac{-4}{\cdot} \right) \end{aligned}$} & 3 & (3, 12), (3, 48) & $\left( \frac{-3}{\cdot} \right)$, $\left( \frac{-3}{\cdot} \right)$\\
 \hline
$ \eta(z) \eta(2z) \eta(7z) \eta(14z)$ & 3 & (2, 14) & $1_{14}$\\
 \hline
\makecell{$\begin{aligned}\eta(4z)^{6}, \:&\frac{\eta(8z)^{18}}{\eta(4z)^6 \eta(16z)^6} \\ &= \eta(4z)^{6} \otimes \left( \frac{-8}{\cdot}\right)\end{aligned}$} & 3 & (3, 16), (3, 64) & $\left( \frac{-4}{\cdot} \right), \left( \frac{-4}{\cdot} \right)$\\
  \hline
$\eta(2z)^{2} \eta(10z)^2$ & 3 & (2, 20) & $1_{20}$\\
   \hline
\makecell{$\begin{aligned}\eta(6z)^{4}, \: &\frac{\eta(12z)^{12}}{\eta(6z)^4 \eta(24z)^4} \\ &= \eta(6z)^{4} \otimes \left( \frac{-4}{\cdot} \right)\end{aligned}$} & 3 & (2, 36), (2, 144) & $1_{6}, 1_{12}$\\
 \hline
\end{tabular}
 % }
  \caption{Parameter Values in Theorem \ref{type2_cong}}
\end{table}

\textbf{Remarks:} 
\begin{enumerate}
\item When $\ell \mid N$, these congruences are equivalent to congruences for coefficients of a form on $\ell$-free level as in the third remark following Theorem \ref{type1_cong}.
\item Let $K$ be an imaginary quadratic field with discriminant $D_k$. We recall that a form $f$ has CM by $\chi = \left( \frac{D_{K}}{\cdot} \right)$ associated to $K$ if and only if $f \otimes \chi = f$. Hence, a CM form has coefficients $a(p) = 0$ when $p$ is inert in $K$, a set of primes of Dirichlet density 1/2. When $f$ is a CM newform, Serre's work in \cite{serre1985lacunarite} implies that almost all of its coefficients are zero. The eta-quotient $\eta(z)^3 \eta(7z)^3, \eta(3z)^8$, and the three eta-quotients $\eta(2z)^3 \eta(6z)^3$, $\eta(4z)^6$, $\eta(6z)^4$ and their twists are CM newforms. When $\ell$ is exceptional of Type II for one of these forms, its coefficients vanish modulo $\ell$ for a set of Dirichlet density $>$ 1/2. We further remark that if $f$ is a newform with odd weight and  real coefficients, then Corollary 1.2 of \cite{schutt2009cm} asserts that $f$ has quadratic nebentypus character $\psi$ and it has CM by this character. 

\medskip

\item Suppose that $\ell$ is exceptional of Type I with $m^{'} - m = \frac{\ell - 1}{2}$. Then we have, for all primes $p \nmid N \ell$,
\[
    a(p) \equiv p^m (1 + p^{m' - m}) \equiv p^m \left(  1 + p^{\frac{\ell - 1}{2}}\right) \equiv p^m \left( 1 + \left( \frac{p}{\ell} \right) \right) \equiv \begin{cases}
        2 p^m, \qquad \left( \frac{p}{\ell} \right) = 1 \\
        0, \quad \quad \quad \left( \frac{p}{\ell} \right) = -1
    \end{cases} \pmod{\ell}.
\]
Hence, $\ell$ is exceptional of Type II. For example, for $f \in \{\eta(z)^8 \eta(2z)^8, \eta(z)^6 \eta(3z)^6\}$ with $(m, m') = (0,1)$, the prime $\ell = 3$ is exceptional of both Type I and Type II.
\end{enumerate}
\subsection{Extension to prime power modulus}
In \cite{swinnerton1973ℓ}, Swinnerton-Dyer extends congruences to prime power modulus. As an example, for $\Delta(z) = \displaystyle{\sum_{n \geq 1}} \tau(n) q^n \in S_{12}(\Gamma_{0}(1))$, we have 
\begin{align}
\tau(n) &\equiv 1537 \: \sigma_{11}(n) \pmod{2^{12}} \:\: \text{if} \:\: n \equiv 5 \pmod{8} \\
\tau(n) &\equiv n^{-30} \: \sigma_{71}(n) \pmod{5^3} \:\: \text{if} \:\: \gcd(n, 5) = 1. \\
\tau(n) &\equiv n^{-610} \: \sigma_{1231}(n) \begin{cases} \pmod{3^6} \:\: \text{if} \:\: n \equiv 1 \pmod{3}, \\ \pmod{3^7} \:\: \text{if} \:\: n \equiv 2 \pmod{3}.  \end{cases}
\end{align}
Our next result is a level $N > 1$ analogue of these congruences in which we extend our classification of primes exceptional of Type I for eta-quotient newforms $f = \displaystyle{\sum_{n}} a(n) q^n$ in Theorem \ref{type1_cong} to prime power modulus. 

\begin{theorem}\label{type1_primepower}
  Let $\ell$ be an exceptional prime of Type I for an eta-quotient newform $f$. Let $t > 1$ and let $0 \leq m < m^{'} \leq \phi(\ell^t)$ with $m + m^{'} \equiv k - 1 \pmod{\phi(\ell^t)}$. The following table gives congruences of the form,
\begin{equation}\label{prime_powers_type1}
a(p) \equiv p^m + p^{m'} \pmod{\ell^t} \:\: \text{for all} \:\: p \equiv b \pmod{d}.
\end{equation}
\end{theorem} 
  \renewcommand{\arraystretch}{1.0} % Adjust the vertical spacing
\setlength{\tabcolsep}{12pt} % Adjust the horizontal spacing
\begin{table}[H]
  \centering
  %\resizebox{1.3\textwidth}{!}{%
 % \resizebox{\textwidth}{!}{%
\begin{tabular}{|c|c|c|c|c|c|}
\hline
$f(z)$ & $\ell$ & $(m, m')$ & $t$ & $b$ & $d$\\ 
\hline
\centering{$\eta(z)^8 \eta(2z)^8$} & 2 & (0,7) & 6 & \text{all residue classes} & 64 \\
\centering{} & 3 & (12,13) & 3 & \text{all residue classes} & 27 \\
\hline
\centering{$\eta(z)^6 \eta(3z)^6$} & 2 & (0,5) & 4 & 5, 7, 11, 19 & 24 \\
\centering{} & 2 & (0,5) & 5 & 13, 17, 23 & 24 \\
\centering{} & 2 & (0,5) & 6 & 1 & 24 \\
\hline
\centering{$\eta(2z)^{12}$} & 2 & (0,5) & 8 & 3 & 8 \\
\centering{} & 2 & (0,5) & 9 & 7 & 8 \\
\centering{} & 2 & (0,5) & 10 & 5 & 8 \\
\centering{} & 2 & (0,5) & 11 & 1 & 8 \\
\centering{} & 3 & (1,4) & 2 & \text{2, 5} & 9 \\
& 3 & (1, 4) & 3 & 8, 17, 26  & 27 \\
\hline
\centering{$\eta(z)^{4} \eta(5z)^{4}$} & 5 & (1,2) & 2 & 1,6,7,11,16,18,21,24 & 25 \\
\hline
\centering{$\eta(3z)^{8}$} & 2 & (0,1) & 2 & 3 & 4 \\
\centering{} & 3 & (0,3) & 4 & \makecell{1, 4, 7, 10, 13, 16, 19, 22, \\ 25,26,28, 31, 34, 37, 40, 43, \\  46, 49, 52, 53, 55, 58, 61, 64, \\ 67, 70, 73, 76, 79, 80}  & 81 \\
\hline
\centering{$\eta(z)^{2} \eta(2z)^{2} \eta(3z)^2 \eta(6z)^2$} & 2 & (0,1) & 2 & \text{all residue classes} & 4 \\
\hline
\centering{$\eta(z)^{2} \eta(11z)^{2}$} & 5 & (0,1) & 2 & 1, 6, 11, 16, 21 & 25 \\
\hline
\centering{$\eta(z) \eta(2z) \eta(7z) \eta(14z)$} & 3 & (0,1) & 2 & 1, 4, 7 & 9 \\
\hline
\centering{$\eta(z) \eta(3z) \eta(5z) \eta(15z)$} & 2 & (0,1) & 3 & \text{all residue classes} & 8 \\
\hline
\centering{$\eta(3z)^2 \eta(9z)^2$} & 3 & (0,1) & 3 & 1, 10, 19, 26 & 27 \\
\hline
\end{tabular}
 % }
 \caption{Parameter Values in \eqref{prime_powers_type1}} \label{type1_params}
\end{table}
%\newpage
\noindent
\textbf{Remarks:}
We also obtain congruences of the form $a(p) \equiv c \: (p^m + p^{m'}) \pmod{e}$ as we illustrate with the following examples:
\begin{enumerate}
    \item Let $f(z) = \eta(2z)^{12}$. For every odd prime $p$, there exists a 2-adic unit $u$ depending on the residue class of $p$ modulo 8, such that
    \[
    a(p) \equiv u (1 + p^5) \pmod{2^{11}}.
    \]
    In particular, we have
    \[
    u =  \begin{cases} 1, \quad \quad \;\; p \equiv 1 \pmod{8} \\
    1729, \quad p \equiv 3 \pmod{8} \\
       1537, \quad  p \equiv 5 \pmod{8} \\
       193, \quad \;\;  p \equiv 7 \pmod{8}
    \end{cases} \quad \text{with modulus sharpening to} \quad \begin{cases} 2^{12}, \quad p \equiv 3,5 \pmod{8} \\
    2^{14}, \qquad p \equiv 7 \pmod{8}.
    \end{cases}
    \]

    \item Let $f(z) = \eta(z)^{6} \eta(3z)^6$. Then we have
    \[
    a(p) \equiv \begin{cases} 
    5 (1 + p^5) \pmod{2^5}, \quad p \equiv 11,19 \\
    9 (1 + p^5) \pmod{2^5}, \quad p \equiv 5     
    \end{cases} \pmod{24}.
    \]
\end{enumerate}

Our next result focuses on extending Theorem \ref{type2_cong} to prime power modulus, just as we extended Theorem \ref{type1_cong}.
\begin{theorem} \label{type2_primepower}
Let $\ell$ be an exceptional prime of Type II for an eta-quotient newform $f$. Let $a$ be a positive integer. The following table gives all congruences of the form
      \[
    f \otimes 1_{\ell} \equiv f \otimes \left( \frac{\cdot}{\ell} \right) \pmod{\ell^{a}}.
      \]
      \end{theorem}
  \renewcommand{\arraystretch}{1.5} % Adjust the vertical spacing
\setlength{\tabcolsep}{12pt} % Adjust the horizontal spacing
\begin{table}[ht]
  \centering
  \small
  %\resizebox{1.3\textwidth}{!}{%
 % \resizebox{\textwidth}{!}{%
\begin{tabular}{|c|c|c|}
\hline
$f(z)$ & $\ell$ & $a$ \\ 
\hline
\centering{$\eta(2z)^{12}$} & 3 & 1,2,3 \\
\hline
\centering{$\eta(3z)^{8}$} & 3 & $\geq 1$ \\
\hline
\centering{$\eta(2z)^{3} \eta(6z)^3$} & 3 & $\geq 1$ \\
\hline
\centering{$\eta(3z)^{2} \eta(9z)^2$} & 3 & $\geq 1$ \\
\hline
\centering{$\eta(6z)^{4}$} & 3 & $\geq 1$ \\
\hline
\end{tabular}
 % }
 % \caption{Parameter Values in Theorem \ref{type2_cong} for % $f$ and $\ell$}
\end{table}
\subsection{Plan for the paper}
The outline for the rest of the paper is as follows. In Section \ref{background}, we give the necessary background, and in Section 3, we prove Theorem \ref{type1_cong}, Theorem \ref{type2_cong}, Theorem \ref{type1_primepower}, and Theorem \ref{type2_primepower}.

\section{Acknowledgment}
The authors are grateful to Matthew Boylan for suggesting this project and for many insightful discussions and comments. We thank Jane Street and the University of South Carolina for their support.

\section{Background} \label{background}
For background on modular forms modulo $\ell$ and level $N = 1$, we refer the reader to \cite{swinnerton1973ℓ} and for levels $N \geq 4$, we refer the reader to \cite{10.1007/BFb0063944, 10.1215/S0012-7094-90-06119-8}.  For details on classical modular forms, one may consult \cite{book}.
Let $\ell$ be prime, let $\Z_{(\ell)}$ denote the localization of $\Z$ at $\ell$, let $N\geq 1$, and let $\chi$ be a Dirichlet character modulo $N$. For $k \geq 0$, let
\begin{equation}\label{l-integrality}
M_k(\Gamma_0(N), \chi)_{(\ell)} = 
M_k(\Gamma_0(N), \chi) \cap \Z_{(\ell)}\llbracket q \rrbracket
\end{equation}
denote the space of weight $k$ modular forms on $\Gamma_0(N)$ with
coefficients in $\Z_{(\ell)}$. If $f(z) = \sum a_n q^n \in \Z_{(\ell)}\llbracket q \rrbracket$ then we define

\[
\widetilde{f} = \sum \widetilde{a_n}\, q^n \in \F_{\ell}[[q]]
\]
to be its coefficient-wise reduction modulo $\ell$, and we define 
\[
\widetilde M^{(\ell)}(\Gamma_0(N)) =
\{\widetilde{f} : f \in M_k(\Gamma_0(N), \chi)_{(\ell)}\}
\subseteq \F_{\ell}\llbracket q\rrbracket,
\]
the algebra of holomorphic integer weight modular forms on $\Gamma_0(N)$ with coefficients in $\Q_{(\ell)}$ reduced modulo $\ell$. 
For $\widetilde{f} \neq 0$, the weight filtration is
\begin{equation}
    w_{\ell}(\widetilde{f}) = \min\{k^{\prime} \geq 0: \:\: \text{there exists} \:\: g \in M_{k^{\prime}}(\Gamma_{0}(N), \chi)_{(\ell)} \:\: \text{with} \:\: \widetilde{f} = \widetilde{g} \}.
\end{equation}
For $\widetilde{f} = 0$, we set $w_{\ell}(\widetilde{f}) = -\infty$. Moreover, if $0 \neq f \in M_{k}(\Gamma_{0}(N), \chi)_{(\ell)}$, then $w_{\ell}(\tilde{f}) \equiv k \pmod{\ell - 1}$ since $\widetilde{E}_{\ell - 1} = 1$. For $t > 1$, the weight filtration modulo $\ell^t$ is well-defined up to $\ell^{t - 1} (\ell - 1)$ \cite{serre1973formes}.

\medskip
Now, we state an important theorem of Sturm \cite{10.1007/BFb0072985} which provides a method to test whether two modular forms with $\ell$-integral coefficients are congruent modulo a prime $\ell$. Let $\ord_{\ell}(f(z)) = \min\{ n \geq 0: a(n) \neq 0 \pmod{\ell}\}$. Buzzard \cite{buzzard2002mod} showed that Sturm Bound for forms with non-trivial Nebentypus agrees with Sturm bound for $\Gamma_{0}(N)$.
 
\begin{theorem}[Sturm] \label{SB}
     Let $N \geq 1$, let $\ell$ be prime and let $f(z), g(z) \in M_{k}(\Gamma_{0}(N), \chi)_{(\ell)}$. Assume that
     \[
     \textup{ord}_{\ell}(f - g) > \begin{cases}  \frac{kb}{12}  - \frac{b - 1}{N}, \quad f - g \in S_{k}(\Gamma_{0}(N), \chi)_{(\ell)} \\ \frac{kb}{12}, \qquad \qquad \text{otherwise}
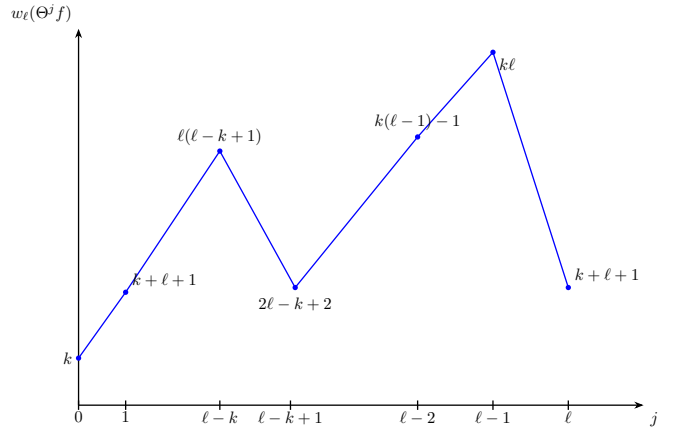
     \end{cases} \: \text{where} \;\:\: b = [\Gamma_{0}(1) : \Gamma_{0}(N)] = N \prod_{p \: \text{prime}\:: \: p \mid N} \left(1 + \frac{1}{p}\right),
     \]
Then we have $\textup{ord}_{\ell}(f(z) - g(z)) = \infty$.
\end{theorem}
\noindent
The same bound works modulo $\ell^t$ for $t > 1$ as shown in \S3.3.~of \cite{taixes2010computing}.
\medskip

We also require generalized Eisenstein Series.
Let $k \geq 4$ and $N \geq 1$. For $\psi$ and $\phi$ modulo $N$, we define the generalized divisor sum as follows:
\begin{equation}
 \sigma_{k}^{\psi, \phi}(n) = \sum_{d\mid n} \psi\left(\frac{n}{d}\right) \phi(d) d^{k}.
\end{equation}
\noindent
If at least one of $\psi, \phi$ is primitive, and $\psi(-1) \phi(-1) = (-1)^k$, then we define the generalized Eisenstein series by: 
\begin{equation}\label{gen_eisenstein}
      G^{\psi, \phi}_k(z) = \begin{cases}
        -\frac{B_{k,\phi}}{2k} + \displaystyle{\sum_{n=1}^\infty} \sigma_{k-1}^{1_1, \phi}(n) q^n, & \text{ if } \psi = 1_1 \\
        \displaystyle{\sum_{n=1}^\infty} \sigma_{k-1}^{\psi, \phi}(n) q^n, & \text{ otherwise};
    \end{cases}
\end{equation}
where $B_{k, \phi}$ denotes the generalized Bernoulli numbers attached to $\phi$. Moreover, we have $G^{\psi, \phi}_k(z) \in M_{k}(\Gamma_{0}(N^2), \psi \phi)$. We call  $G^{\psi, \phi}_k(z)$ to be normalized if $a(1) = 1$, and define $E_{k}^{\psi, \phi} = \frac{G_{k}^{\psi, \phi}}{-B_{k, \phi}/2k}$. 
When $\psi = \phi = 1_1$, then we have $B_{k, \phi} = B_k$, $G_{k}^{\psi, \phi} = G_k$ and $E_{k}^{\psi, \phi} = E_k$. Let $k = 2$ and $N \geq 2$.  Let
\[
E_{2}(z) =  1 - 24 \sum_{n \geq 1} \sigma_{1}(n) q^n
\]
denote the quasimodular Eisenstein series of weight $k = 2$ on $\SL_{2}(\Z)$. The standard weight 2 Eisenstein series of level $N$ is defined as
\begin{equation}\label{defn_E2N}
E_{2, N }(z) = \frac{1}{24}(N E_{2}(Nz) - E_{2}(z)) \in M_{2}(\Gamma_{0}(N)).   
\end{equation}
\noindent
We shift our attention to operators acting on spaces of modular forms. 

Let $\displaystyle{g(z) = \sum_{n \geq 0} b(n) q^n} \in M_{k}(\Gamma_{0}(N), \chi)$. For all integers $m\geq 1$, we define the $U$- and $V$-operators by 
\begin{equation} \label{UVq}
g\mid U_m = \sum_n a(mn)q^n, \ \ \ 
g\mid V_m = \sum_n a(n)q^{mn}. 
\end{equation}
These operators map spaces of modular forms as follows: 
\begin{align}
& V_m : M_k(\Gamma_{0}(N),\chi) \longrightarrow M_k(\Gamma_{0}(mN), \chi), \label{vmap} \\
& U_m : M_k(\Gamma_{0}(N),\chi)\longrightarrow \begin{cases} M_k(\Gamma_{0}(mN), \chi), & m\nmid N, \\ M_k(\Gamma_{0}(N), \chi), & m\mid N.\notag
\end{cases}
\end{align}

We will use the following elementary property of the $U$-operator on $q$-series.
\begin{proposition}\label{U_factor}
Let $m\geq 1$ in $\mathbb{Z}$, and let $f(q)$, $g(q)\in \mathbb{C}[[q]]$.  Then we have 
\[
(f(q)g(q^m))\mid U_{m} = f(q)\mid U_{m} \cdot g(q).
\]  
\end{proposition}

\medskip
\noindent
We now introduce Hecke operators. For $g(z) = \displaystyle{\sum_{n \geq 0}} b(n) q^{n} \in M_k(\Gamma_{0}(N),\chi)$ and a prime $p$, we define the action of Hecke operator $T_p$ on $g(z)$~by
\begin{equation}\label{defn_Hecke}
g \mid T_{p} = g \mid U_p + \chi(p) p^{k - 1} g \mid V_p = \sum_{n \geq 0} \left(b(p n) + \chi(p){p}^{k - 1}b\left(\frac{n}{p}\right)\right)q^n,
\end{equation}
where $b\left(\frac{n}{p}\right) = 0$ when $p\nmid n$. In general, for a positive integer $n$, we define the action of $T_{n}$ on $f(z)$~by
\begin{equation}\label{defn_Tn}
g \mid T_{n} = \sum_{n \geq 0} \left( \sum_{d \mid (m, n)} \chi(d) d^{k - 1} b\left(\frac{mn}{d^2}\right) \right) q^{n}.
\end{equation}
The Hecke operators preserve both $M_k(\Gamma_{0}(N),\chi)$ and its subspaces of Eisenstein series and cusp forms. 

\medskip

We say that $g \in M_{k}(\Gamma_{0}(N),\chi)$ is a Hecke eigenform if for every integer $n \geq 1$ with $\gcd(n, N) = 1$, there exists $\lambda_n \in \C$ with $g \mid_k T_n = \lambda_n g$. Further, $f$ is said to be normalized if $b(1) = 1$.
A newform in $S_{k}^{\new}(\Gamma_{0}(N), \chi)$ is a normalized cusp form that is an eigenform of all the Hecke operators $T_n$ and all of the Atkin-Lehner operators $W_{p}^{N}$ for primes $p \mid N$ and $H_N$, the Fricke involution.

\medskip

The Fourier coefficients of normalized Hecke eigenforms satisfy the following multiplicative relations. 
\begin{proposition}
    Let $\displaystyle{g(z) = \sum_{n} b(n)q^{n}} \in M_{k}(\Gamma_{0}(N), \chi)$ be a normalized Hecke eigenform. Then for all $m, n \geq 0$, we have
    \begin{equation}\label{multiplicativity}
    b(m)b(n) = \sum_{d \mid \gcd(m , n)} \chi(d) d^{k - 1} \left(\frac{mn}{d^2}\right).
    \end{equation}
\end{proposition}    
    \par
    
    We now proceed to the action of Ramanujan's differential operator. The Ramanujan $\theta$-operator is defined by: 
    \begin{equation}\label{defn_theta}
    \theta = q \frac{d}{dq} = \frac{1}{2 \pi i} \frac{d}{dz} \: \text{where} \: q = e^{2 \pi i z}, \: \text{and} \:\:
    \displaystyle{\theta(f(z)) = \theta \left(\sum_{n \geq 0} a(n)q^{n} \right) = \sum_{n \geq 0} n a(n) q^{n}}.
    \end{equation}

The derivative operator fails to preserve modularity in general. However, for a prime $\ell > 3$ and $N = 1$, Lemma 5(ii) of \cite{swinnerton1973ℓ} implies that
\begin{equation}
\theta = q\frac{d}{dq} \, : \,\tilde{M}_{k}^{(\ell)}(\Gamma_0(1)) \longrightarrow \tilde{M}_{k + \ell + 1}^{(\ell)}(\Gamma_{0}(1)). \label{theta_operator}
\end{equation}
Moreover, we have 
\begin{equation}\label{theta_modl}
w_{\ell}(\theta f) = w_{\ell}(f) + \ell + 1 - \alpha (\ell - 1) \; \text{with} \; \alpha \geq 0, \; \text{and} \; \alpha = 0 \; \text{when} \; \ell \nmid  w_{\ell}(f).
\end{equation}
Let $t > 1$. Theorem 1.1 of Chen and Kiming \cite{Chen_Kiming_2016} implies that 
\begin{equation}\label{chen_kiming}
\theta: M_{k}(\Gamma_{0}(N)) \rightarrow M_{k + k(t)}(\Gamma_{0}(N)) \pmod{\ell^t} \:\: \text{with} \:\:  k(t) = 2 + 2 \ell^{t - 1} (\ell - 1),
\end{equation}
which is useful for proving congruences modulo prime powers.
Further, the results of \cite{swinnerton1973ℓ} have been generalized by Katz \cite{10.1007/BFb0063944} and Gross \cite{10.1215/S0012-7094-90-06119-8} to forms of level $N \geq 4$, for $k \geq 2$ and for primes $\ell$ with $\ell \nmid N$.

\section{Proofs of Theorem \ref{type1_cong} and Theorem \ref{type2_cong}}
\subsection{Proof of Type I Congruences} \label{proof_Type1}
Let $N \geq 1$, and let $k \geq 2$ be even. This section is divided into two parts. We first consider the case when $\ell \geq 5$ is prime and $N \geq 4$. We adapt and extend the arguments of \cite{swinnerton1973ℓ, Ribet1975, ribet1985adic}. The central idea involves the use of filtration arguments to bound the possible $\ell$ that could be exceptional. In the later part, we explicitly prove congruences mod $\ell$, and modulo $\ell^t$ for $t > 1$, $\ell \in \{ 2,3 \}$, and $N = \{1,2,3\}$. 

\medskip

\noindent
Let $\psi, \phi$ be real-valued Dirichlet characters mod $N$, and let $\chi = 1_{N}$. Let $\ell \geq 3$ be prime, and let $\gcd(\ell, N) = 1$. We note that for any $a \in (\mathbb{Z}/\ell\mathbb{Z})^\times$ and $b \in (\mathbb{Z}/N\mathbb{Z})^\times$,  there exist infinitely many primes $p$ with $p \equiv a \pmod \ell$, and $ p \equiv b  \pmod N$.
Further, using \eqref{Detcond}, it follows that $m + m' \equiv k - 1 \pmod{\ell - 1}$, $\psi \phi = 1_{N}$, and $m \neq m'$.
\iffalse
\begin{proof}
Let $a$ be a primitive root modulo $\ell$. Using Lemma \ref{inf_p}, there exists a prime $p$ with $p \equiv a \pmod{\ell}$ and $p \equiv 1 \pmod{N}$. In particular, $p \nmid N \ell$, therefore by \eqref{Detcond}, we have $p^{m ' + m} \equiv p^{k - 1} \pmod{\ell}$. Since $p \equiv a \pmod{\ell}$ and $a$ is primitive, we obtain the desired result.
\end{proof}
Our next lemma shows that $\psi$ and $\phi$ are inverses of each other.
\begin{lemma}\label{trivial_char}
 Let $\ell \geq 3$ be prime. Then we have $\psi \phi = 1_{N}$.
\end{lemma}
\begin{proof}
Fix any $b \in (\mathbb{Z}/N\mathbb{Z})^\times$. Using  Lemma \ref{inf_p} and Lemma \ref{mm'res}, we deduce that $\psi(b) \phi(b) \equiv 1 \pmod{\ell}$. Since $\psi(b), \phi(b) \in \pm{1}$, and $\ell \geq 3$, this implies that $\psi(b) \phi(b) = 1$. Since $b$ is arbitrary, we obtain the desired result.
\end{proof}
\begin{lemma}
If $k$ is even, then $m \neq m'$. In particular, we have $2 \leq m' - m + 1 \leq \ell - 1$.
\end{lemma}
\begin{proof}
 If $m = m'$, then Lemma \ref{mm'res} gives $k - 1 \equiv 2m \pmod{\ell - 1}$. Since $\ell \geq 3$, this implies that $k$ is odd, which gives a contradiction. The bounds are immediate from $0 \leq m , m' \leq \ell - 2$.   
\end{proof}
\fi
\noindent
We now proceed to the proof of Theorem \ref{type1_cong}. 

\medskip

We suppose that the eta-quotient newform $f(z) = \sum_{n} a(n) q^n \in S_{k}(\Gamma_{0}(N), \chi)$ satisfies a type I congruence mod $\ell$. We prove part (1) of Theorem \ref{type1_cong} and the proof of part (2) follows similarly by using $G_{2} \equiv G_{\ell + 1} \pmod{\ell}$. Let $3 \leq m' - m + 1 \leq \ell - 2$. Let $E = G_{m' - m + 1} \otimes \psi1_N$. We note that for all primes $p \nmid N \ell$, using \eqref{Trcond}, we have
\begin{equation}\label{eq1}
    a(p) \equiv p^m \: \psi(p)( 1 + p^{m' - m}) \equiv p^m \: \psi(p) \: \sigma_{m' - m}(p)  \pmod{\ell}.
\end{equation}
Multiplying both sides of \eqref{eq1} by $p 1_{N}(p)$, we have
\begin{equation}
    p \: 1_{N}(p) \: a(p) \equiv p^{m + 1} \: 1_{N}(p) \: \psi(p) \: \sigma_{m' - m}(p) \pmod{\ell}.
\end{equation}
for all primes $p \nmid N$. It follows that the coefficients of $q^p$ in $\theta(f \otimes 1_N)$ and $\theta^{m + 1}(E)$ agree modulo $\ell$. Since $f$ and $E$ are normalized Hecke eigenforms for all $T_p$, using \eqref{multiplicativity} and applying induction on $n$ gives the following
relation for all $n \geq 1$:
\begin{equation}\label{eq6}
     n \: 1_{N}(n) \: a(n) \equiv n^{m + 1} \: 1_{N}(n) \: \psi(n) \: \sigma_{m' - m}(n) \pmod{\ell}.
\end{equation}
\par
Thus, using \eqref{theta_operator} and \eqref{eq6}, the coefficient of $q^n$ of $\theta(f \otimes 1_N)$ and $\theta^{m + 1} (G_{m' - m + 1} \otimes \psi1_N)$ agree modulo $\ell$ for each $n \geq 1$. Hence, we conclude that
\begin{equation}\label{cong1_form}
    \theta(f \otimes 1_N) \equiv \theta^{m + 1} (G_{m' - m + 1} \otimes \psi1_N) \pmod{\ell}.
\end{equation}
To show that $\ell < k$, we first claim that $m' + m \ell + 1 \leq k$. Using induction on $j$, it follows that for all $0 \leq j \leq m$, we have $w_{\ell}(\theta^j E) = w_{\ell}(E) + j(\ell + 1) \not\equiv 0 \pmod{\ell}$. Therefore, we have
\begin{equation}\label{eq:LHS}
w_{\ell}(\theta^{m + 1}E) = w_{\ell}(\theta^m E) + \ell + 1 = m' - m + 1 + (m + 1)(\ell + 1).
\end{equation}
On the other hand, we have
\begin{equation}\label{eq:RHS}
w_{\ell}(\theta(f \otimes 1_N)) \leq w_{\ell}(\theta f) \leq k + \ell + 1. 
\end{equation}
Comparing \eqref{eq:LHS} and \eqref{eq:RHS} gives $m' + m\ell + 1 \leq k$.
If $m \neq 0$ then $m' > m \geq 1$ implies that $\ell \leq k - 3 < k$, as desired.
\par We next claim that for $m = 0$, either $\ell < k$ or $\ell \mid (a(p) - \psi(p) \sigma_{k - 1}(p))$. Suppose that $\ell > k$. Since $m = 0$ and $m' + m \ell + 1 \leq k$, we have $m' \leq k - 1$. We claim that $m' = k - 1$. To see this, we suppose that $m' < k - 1$. Since $m = 0$, we have $m' \equiv k - 1 \pmod{\ell - 1}$, so there exists $t \in \Z$ such that $k - 1 - m' = t (\ell - 1)$. Since $m' < k - 1$, we have $t \geq 1$ and $k - 1 - m' \geq \ell - 1$. Since $\ell > k$, we have $\ell - 1 > k - 1 > k - 1 - m'$, which is a contradiction. Hence, we conclude that $m' = k - 1$. Using \eqref{cong1_form}, we deduce that 
\begin{equation}\label{m0_cong}
\theta(f \otimes 1_N) \equiv \theta (G_{k} \otimes \psi 1_N) \pmod{\ell}. 
\end{equation}
\par
We next want to show that $f \otimes 1_N \equiv (G_k \otimes \psi) \otimes 1_N \pmod \ell.$ We argue by contradiction. Since $f \otimes 1_N  - ((G_k \otimes \psi) \otimes 1_N) \in M_k$ and $f \otimes 1_N  - ((G_k \otimes \psi) \otimes 1_N) \not\equiv 0 \pmod \ell,$ we have
        \[
        0 \leq w_\ell(f \otimes 1_N  - (G_k \otimes \psi 1_N)) \leq k \leq \ell - 2 < \ell - 1,
        \]
        and
        \[
        w_\ell(f \otimes 1_N  - (G_k \otimes \psi 1_N)) \equiv k \pmod{\ell - 1}.
        \]
        Thus, it follows that $w_\ell (f \otimes 1_N  - (G_k \otimes \psi 1_N)) = k \not\equiv 0 \pmod{\ell}$,
        which gives us $w_\ell(\theta(f \otimes 1_N  - (G_k \otimes \psi 1_N))) = k + \ell + 1\neq -\infty$. Hence, we find that
        \[
        \theta(f \otimes 1_N) \not\equiv \theta (G_k \otimes \psi 1_N) \pmod \ell,
        \]
        which contradicts \eqref{m0_cong}.

Therefore, comparing the coefficient of $q^n$ on both sides of $f \otimes 1_N \equiv (G_k \otimes \psi) \otimes 1_N \pmod \ell$  gives
\begin{equation}\label{eq4}
    1_{N}(n) \: a(n) \equiv \psi(n) \: 1_{N}(n) \: \sigma_{k - 1}(n) \pmod{\ell}. 
\end{equation}
Let $\gcd(p, N) = 1$. Substituting $n = p$ in \eqref{eq4}, we obtain $a(p) \equiv \psi(p) \: \sigma_{k - 1}(p) \pmod{\ell}$.

\medskip

We now address the case when $\ell \in \{2,3\}$ is exceptional of Type I and $N \in \{1, 2,3\}$. To prove such congruences modulo $\ell$, we use elementary $q$-series identities. We prove the result for $f(z) = \eta(z)^8 \eta(2z)^8 \\ \in S_{8}(\Gamma_0(2))$. The proofs in remaining cases follow similarly. Let
\begin{equation}\label{defn_G2}
 G_2(z) = E_{2, 2}(z) = \frac{1}{24}(2 E_{2}(2z) - E_{2}(z)) = \frac{1}{24} + \sum_{n \geq 1} \left(\sigma_{1}(n) - 2 \sigma_{1}\left( \frac{n}{2} \right)  \right) q^n \:\: \text{where} \:\: \sigma_{1}\left( \frac{n}{2} \right) = 0 \:\: \text{for} \:\: 2 \nmid n. 
\end{equation}We aim to show, for $\ell \in \{2,3\}$ that 
\begin{equation}\label{type1_mod2}
\theta(f \otimes 1_2) \equiv \theta(G_{2}(z) \otimes 1_2) \pmod{\ell}.
\end{equation}
For $\ell = 2$, \eqref{type1_mod2} reduces
\begin{equation*}
    \theta f(z) \equiv \theta G_{2}(z) \pmod{2} 
\end{equation*}
We note using Jacobi Triple Product Identity that
\begin{equation}
\Delta(z) \equiv \eta(8z)^3 = \sum_{n \geq 0}(-1)^n (2n + 1)q^{(2n + 1)^2} \equiv \sum_{n \geq 0} q^{(2n + 1)^2} \pmod{2}.
\end{equation}
Using $f(z) \equiv \Delta(z) \pmod{2}$,
\begin{equation}\label{theta_f}
\theta f(z) \equiv \theta\Delta(z) \equiv \sum_{n \geq 0} (2n + 1)^2 q^{(2n + 1)^2} \equiv \sum_{n \geq 0} q^{(2n + 1)^2} \pmod{2}.
\end{equation}
On the other hand, using \eqref{defn_G2}, we have
\begin{equation}\label{theta_G2}
\theta G_{2}(z) \equiv \sum_{n \geq 0} n \sigma_{1}(n) q^n = \sum_{n \geq 0} q^{(2n + 1)^2} \pmod{2}
\end{equation}
The last equality follows since the sum in \eqref{theta_G2} vanishes for even $n$, and moreover, $\displaystyle{\sigma_{1}(n) = \prod_{p^a || n}} (1 + a) \pmod{2}$ is odd if and only if $a$ is even. Comparing \eqref{theta_f} and \eqref{theta_G2}, we conclude that $\theta f(z) \equiv \theta G_{2}(z) \pmod{2}$.

\medskip

We turn to the proof of $\theta (f(z) \otimes 1_2) \equiv \theta (G_{2}(z) \otimes 1_{2}) \pmod{3}$. Using \eqref{defn_G2}, it follows that 
\begin{equation}\label{LHS}
    \theta(G_{2}(z) \otimes 1_2) = \sum_{\substack{n \geq 1 \\ n \: \text{odd}}} n \sigma_{1}(n) q^n \pmod{3}.
\end{equation}
Let $F(z) = \frac{1}{480} (E_8(z) - 2^7 E_{8}(2z)) \in M_{8}(\Gamma_{0}(2))$. Then we have $f(z) \otimes 1_3 \equiv F(z) \otimes 1_3 \pmod{3}$, and for $f(z) = \displaystyle{\sum_{n}} a(n) q^n$, we deduce that
\begin{equation} \label{eisen_exp}
    a(n) \equiv \sigma_{7}(n) - 2^7 \sigma_{7}\left( \frac{n}{2} \right) \pmod{3} \:\: \text{for all} \:\: (n, 3) = 1.
\end{equation}
By Fermat's Little Theorem, we deduce that $\sigma_{7}(n) \equiv \sigma_{1}(n) \pmod{3}$. Further, using \eqref{eisen_exp}, for odd $n$ with $(n, 3) = 1$, we have $a(n) \equiv \sigma_{1}(n) \pmod{3}$. Hence, we conclude that for all odd $n$, we have
\begin{equation}
  \theta(G_{2}(z) \otimes 1_2) = \sum_{\substack{n \geq 1 \\n \: \text{odd}}} n \sigma_{1}(n) q^n \equiv \sum_{\substack{n \geq 1 \\n \: \text{odd}}} n a(n) q^n \pmod{3} \equiv \theta(f(z) \otimes 1_2) \pmod{3}. 
\end{equation}
Next, to prove a Type I congruence for a newform $f$ modulo prime powers, we require the following theorem that helps in determining the action of theta operator on $f$ modulo powers of 2 and 3.
\begin{theorem}\label{theta_2_3_powers}
Let $N \geq 1$ and let $\ell \in \{2,3\}$. We define 
\begin{equation}
    j = j_{\ell, t} = \begin{cases}
        2 + \phi(2^t), \quad \text{if} \:\: \ell = 2 \:\: \text{and} \:\: t \geq 4 \\
        2 + \phi(3^t), \quad \text{if} \:\: \ell = 3 \:\: \text{and} \:\: t \geq 2. 
    \end{cases}
\end{equation}
Then there exists an Eisenstein series
\[
  F_{\ell, t}(z) \in M_{j}\bigl(\Gamma_0(\ell^{t - 1})\bigr)
\]
such that
\[
  E_2(z) \equiv F_{\ell,t}(z) \pmod{\ell^t}.
\]
Consequently, for every $f\in M_k(\Gamma_0(N), \chi)$ one has the congruence
\[
  \theta f = \vartheta f - \frac{k}{12} E_{2} f \equiv \vartheta f - \frac{k}{12} F_{\ell,t} f \pmod{\ell^t}
\]
Hence, $\theta f \in M_{j + k}(\Gamma_0(\ell^{t - 1} N), \chi) \pmod{\ell^t}$. 
\end{theorem}

\begin{proof}
Let $\ell = 3$, let $t \geq 2$, let $j = 2 + \phi(3^t)$ and let $j^{\prime} = 2$. Since $j \equiv j^{\prime} \pmod{ \phi(3^t)}$, by Kummer's congruence [Proposition 11.4.4, \cite{cohen2007number}], we have 
\begin{equation} \label{Kummer}
(3^{1 + \phi(3^t)} - 1) \: \frac{B_{2 + \phi(3^t)}}{{2 + \phi(3^t)}} + \frac{2/3}{2 + \phi(3^t)} - (2 + \phi(3^t)) \equiv B_{2} + \frac{1}{3} - 2 \pmod{3^t}.
\end{equation}
Further, using the Clausen-von Staudt Theorem [Theorem 9.5.14, \cite{cohen2007number}], we deduce that $v_{3}(B_{2 + \phi(3^t)}) = -1$. Hence, we have 
\begin{equation}\label{3-adic_val}
v_{3}\left( 3^{1 + \phi(3^t)}  \: \frac{B_{2 + \phi(3^t)}}{{2 + \phi(3^t)}} \right) = \phi(3^t) \geq t. 
\end{equation}
Substituting \eqref{3-adic_val} in \eqref{Kummer} and simplifying gives 
\begin{equation}\label{kummer_simp}
- \frac{B_{2 + \phi(3^t)}}{{2 + \phi(3^t)}} \equiv  B_2 + \frac{3^{t - 2}}{1 + 3^{t -1}} + 2 \cdot  3^{t - 1} \pmod{3^t}.
\end{equation}
On multiplying both sides of \eqref{kummer_simp} by 3, we obtain the following:
\begin{align}\label{cong_3^{t - 1}}
  \frac{-3 B_{2 + \phi(3^t)}}{{2 + \phi(3^t)}} \equiv 3 B_{2} \pmod{3^{t - 1}}.
\end{align}
Since the terms on either side of the congruence are units in $\Z_{3}$, taking inverses in \eqref{cong_3^{t - 1}} and multiplying by 6 yields
\begin{equation}
    \frac{ - 2 (2 + \phi(3^t))}{\: B_{2 + \phi(3^t)}} \equiv \frac{2}{B_{2}} \pmod{3^{t}}.
\end{equation}
We recall the $q$‑expansion of the normalized Eisenstein series  $E_{2 + \phi(3^t)}(z) \in M_{2 + \phi(3^t)}(\Gamma_{0}(1))$:
\begin{align}
  E_{2 + \phi(3^t)}(z) &= 1 - \frac{2 (2 + \phi(3^t)) }{B_{2 + \phi(3^t)}} \: \sum_{n \geq 1}\left(\sum_{d \mid n} d^{1 + \phi(3^t)} \right) q^n   \\
  &\equiv 1 - \frac{2 (2 + \phi(3^t)) }{B_{2 + \phi(3^t)}} \: \sum_{n \geq 1}\left(\sum_{\substack{d \mid n \\ 3 \nmid d}} d \right) q^n \pmod{3^t}. \label{eq:k}
\end{align}
Further, using \eqref{gen_eisenstein}, taking $\psi = 1_1$, and $\phi = 1_3$, we find that
\begin{align}
 E_{2}^{1_1, 1_3}(z) = 1 -  \frac{4}{B_{2, 1_3}} \sum_{n \geq 1} \sigma_{1}^{1_1, 1_3}(n) q^n = 1 + \frac{2}{B_2} \sum_{n \geq 1} \left( \sum_{\substack{d \mid n \\ 3 \nmid d}} d \right) q^n \in M_{2}(\Gamma_{0}(3)),
\end{align}
where the last equality follows from $L(s, 1_3) = \zeta(s) (1 - 3^{-s})$ and $\zeta(1 - n) = \frac{-B_n}{n}$ for $n \geq 1$. Moreover, we observe that
\begin{equation}
 E_{2}^{1_1, 1_3}(z) = 1 + \frac{2}{B_2}\sum_{n \geq 1}\left(\sum_{\substack{d \mid n \\ 3 \nmid d}} d \right) q^n
  =  -\frac{1}{2} \left(E_2(z) - 3E_2(3z) \right)
  =  -\frac{1}{2} E_2(z) \big|(1 - 3 V_3) \in M_{2}(\Gamma_{0}(3)). \label{eq:k'}
\end{equation}
Since $j = 2 + \phi(3^t) \equiv 2 = j^{'} \pmod{\phi(3^t)}$, comparing \eqref{eq:k} and \eqref{eq:k'}, we deduce that
\begin{equation}\label{eq:E_k}
    E_{2 + \phi(3^t)}(z) = -\frac{1}{2} E_2(z) \big|(1 - 3 V_3) \pmod{3^t}.
\end{equation}
Acting with the operator
\begin{equation}\label{V3_op}
    (1 - 3 V_3)^{-1} \equiv \sum_{i = 0}^{t - 1} (3V_3)^i \pmod{3^t}
 \end{equation}
 on both sides of \eqref{eq:E_k}, we get
 \begin{equation}\label{simp_cong_3^t}
   E_{2}(z) \equiv -2 \: E_{2 + \phi(3^t)}(z) \bigg| \sum_{i = 0}^{t - 1} (3V_3)^i \equiv -2 \sum_{i = 0}^{t - 1} 3^i \: E_{2 + \phi(3^t)}(3^i z) \pmod{3^t}. 
 \end{equation}
Let $t \geq 2$ and $0 \leq i \leq t - 1$. Let $E_{2}^{*}(z) = -\frac{1}{2} E_2(z) \big|(1 - 3 V_3)$. Then we have
 \begin{align}
     E_{2 + \phi(3^t)}(3^i z) &\equiv E_{2 + \phi(3^{t - i})}(3^i z) \equiv E_{2}^{*}(3^i z) \pmod{3^{t - i}} \\
     3^i \: E_{2 + \phi(3^t)}(3^i z) &\equiv 3^i \: E_{2 + \phi(3^{t - i})}(3^i z) \pmod{3^t}. \label{change_index}
 \end{align}
Using \eqref{change_index}, the equation \eqref{simp_cong_3^t} transforms as
\begin{align}
    E_{2}(z) \equiv - 2 \sum_{i = 0}^{t - 1} 3^i E_{2 + \phi(3^{t - i})} (3^i z) \pmod{3^t}.
\end{align}
Thus, there exists a modular form $F_{3, t}(z) = - 2 \displaystyle{\sum_{i = 0}^{t - 1}} 3^i \: E_{2 + \phi(3^{t})} (3^i z) \in M_{2 + \phi(3^t)}(\Gamma_{0}(3^{t - 1}))$ such that $E_{2}(z) \equiv F_{3,t} \pmod{3^t}$. For $p = 2$, we note that $E_{2} \equiv 1 \pmod{8}$ and for $t \geq 4$, the proof follows similarly by using Kummer's congruence [Proposition 11.4.4, \cite{cohen2007number}]. Let $j = 2 + \phi(2^t) \equiv 2 = j^{'} \pmod{2^{t - 1}}$. Then, we have
\[
(2^{1 + 2^t} - 1) \: \frac{B_{2 + 2^t}}{{2 + 2^t}} + \frac{1/2}{2 + 2^t} \equiv \frac{B_{2}}{2} + \frac{1}{4} \pmod{2^t}.
\]
\end{proof}
We now apply Theorem \ref{theta_2_3_powers} to $f(z) = \eta(z)^8 \eta(2z)^8 \in S_{8}(\Gamma_{0}(2))$ and obtain a Type I congruence modulo $3^3$. The proof in remaining cases mod $\ell^t$ for $\ell \in \{2,3\}$ and $t > 1$ follows similarly.
For $\ell \geq 5$ and $t > 1$, we use \eqref{chen_kiming} and Sturm's bound \eqref{SB} to prove congruences mod $\ell^t$.
\par
Let $f(z) = \sum a(n) q^n$. From Table \ref{type1_params}, we have
\begin{align}
    a(p) &\equiv p^{12} + p^{13} \pmod{27} \:\: \text{for all primes} \:\: p \not\in \{2,3\} \\
    a(p) &\equiv p^{12}(1 + p) \equiv p^{12}(1 + p^{19}) \pmod{27} \; \; \text{using Fermat's Little theorem}. \label{FLT}
\end{align}
Multiplying both sides of \eqref{FLT} by $p^3$, we obtain
\begin{align}
    p^3 a(p) &\equiv p^{15}( 1 + p^{19}) \pmod{27} \;\; \text{for all primes} \;\; p \neq 2 \\
     p^3 a(p) &\otimes 1_2(p) \equiv p^{15}( 1 + p^{19}) \otimes 1_{2}(p) \pmod{27} \;\; \text{for all primes} \;\; p. \label{prime_cong3}
\end{align}
Using \eqref{multiplicativity} and induction on $n$ gives the following relation for all $n \geq 1$:
\begin{equation}
    n^3 a(n) \otimes 1_2(n) \equiv n^{15}( 1 + n^{19}) \otimes 1_{2}(n) \pmod{27}
\end{equation}
Hence, we conclude that
\begin{equation}
    \theta^3(f \otimes 1_2) \equiv \theta^{15}(G_{20} \otimes 1_2) \pmod{27}.
\end{equation}
Using Theorem \ref{theta_2_3_powers}, we obtain that $\theta^3(f \otimes 1_2) \in S_{68}(\Gamma_{0}(36))$ and $\theta^{15}(G_{20} \otimes 1_2) \in S_{320}(\Gamma_{0}(36))$. Since $E_{18} \equiv 1 \pmod{27}$, to prove that
\begin{equation}
     \theta^3(f \otimes 1_2) \: E_{18}^{12} \equiv \theta^{15}(G_{20} \otimes 1_2) \pmod{27}
\end{equation}
in $S_{320}(\Gamma_{0}(36))$, using \ref{SB}, it suffices to check coefficients of index $n$ up to 1918, which we have verified.

\medskip
\noindent
We now proceed to the proof of Theorem \ref{type2_cong}.
\subsection{Proof of Type II Congruences}
Let $\ell > 3$ be prime. Let $\ell$ be an exceptional prime of Type II for $f$. Then for all primes $p \nmid N \ell$ with $\left( \frac{p}{\ell} \right) = -1$, we have $a(p) \equiv 0 \pmod{\ell}$. Consider
\begin{align*}
    p \left[ \left( \frac{p}{\ell} \right) - 1 \right] a(p) \equiv \displaystyle{\begin{cases}
         0 \pmod{\ell} &\text{if $\ell \mid p$}, \\
         0 \pmod{\ell} &\text{if $\left( \frac{p}{\ell} \right) = -1$}, \\
         0 \pmod{\ell}, &\text{if $\left( \frac{p}{\ell} \right) = 1$}.
    \end{cases}}
\end{align*}
In all three cases, we have $p \left[ \left( \frac{p}{\ell} \right) - 1 \right] a(p) \equiv 0 \pmod{\ell}$. Since $\left( \frac{p}{\ell} \right) \equiv p^{\frac{\ell - 1}{2}} \pmod{\ell}$, for all primes 
\begin{equation}
 p^{\frac{\ell + 1}{2}} \cdot 1_{N}(p) \cdot a(p) \equiv p \cdot 1_N(p) \cdot a(p) \pmod{\ell}.
\end{equation}
Since $f$ is a normalized eigenform for all $T_p$, using \eqref{multiplicativity}, induction on $n \geq 0$ gives
\begin{equation}\label{eq7}
 n^{\frac{\ell + 1}{2}} \cdot 1_N(n) \cdot a(n) \equiv n \cdot 1_N(n) \cdot a(n) \pmod{\ell}.    
\end{equation}
Hence, using \eqref{theta_operator} and \eqref{eq7} we conclude that 
\begin{equation}\label{type2_form}
\theta^{\frac{\ell + 1}{2}} (f \otimes 1_N) \equiv \theta (f \otimes 1_N) \pmod{\ell}.
\end{equation}
We next claim that $\ell < 2k$. On the contrary, suppose that $\ell > 2k$. It follows by induction that for $0 \leq j \leq \frac{\ell - 1}{2}$, we have $w_{\ell}(\theta^j (f \otimes 1_N)) = w_{\ell}(\theta^j f) = k + j (\ell + 1) \not\equiv 0 \pmod{\ell}$. 
In particular, we have $w_{\ell}(\theta^{\frac{\ell - 1}{2}}(f \otimes 1_N)) = k + \left(\frac{\ell - 1}{2} \right) (\ell + 1)$ which implies that 
\begin{equation}\label{fil_lhs}
w_{\ell}(\theta^{\frac{\ell + 1}{2}}(f \otimes 1_N)) = k + \left(\frac{\ell + 1}{2} \right)(\ell + 1). 
\end{equation}
We note that $\ell > 2k > k > 3$, implying that $4 \leq k \leq \ell - 2$. Since $f \in S_k(\Gamma_{0}(N), \chi)$ and $w_{\ell}(f \otimes 1_N) = w_{\ell}(f) = k \not\equiv 0 \pmod{\ell}$, it follows that 
\begin{equation}\label{fil_rhs}
w_{\ell}(\theta(f \otimes 1_N)) = k + \ell + 1.
\end{equation}
From \eqref{type2_form}, we see that $w_{\ell}(\theta^{\frac{\ell + 1}{2}} (f \otimes 1_N)) = w_{\ell}(\theta (f \otimes 1_N))$. On comparing \eqref{fil_lhs} and \eqref{fil_rhs}, we obtain $\ell \in \{\pm 1\}$, which contradicts that $\ell > 3$. Hence, we conclude that $\ell < 2k$. 

\medskip

Let $k < \ell < 2k$. To prove that $\ell \in \{2k - 1, 2k - 3\}$, we require some preliminaries. 

\begin{proposition} \label{Ul_filtration}
Let $\ell$ be an exceptional prime of Type II. Then we have either
\begin{enumerate}[label=(\alph*)]
    \item $w_{\ell}(\theta^{\ell - 1} f) = w_{\ell}(f)$ or
    \item $w_{\ell}(\theta^{\ell - 1} f) \equiv 0 \pmod{\ell}$.
\end{enumerate}
Furthermore, (a) holds if and only if $f \mid U_{\ell} \equiv 0 \pmod{\ell}$ and (b) holds if and only if $f \mid U_{\ell} \not\equiv 0 \pmod{\ell}$.

\end{proposition}
\begin{proof}
Let $k < \ell < 2k$. Since $f \in S_{k}(\Gamma_{0}(N), \chi)$, we see that $w_{\ell}(f) = k \not\equiv 0 \pmod{\ell}$. We note that if $w_{\ell}(\theta^{\ell - 1}f) \not\equiv 0 \pmod{\ell}$, then using $\theta^{\ell} f \equiv \theta f \pmod{\ell}$ and \eqref{theta_modl}, we deduce that $w_{\ell}(\theta^{\ell - 1}f) = w_{\ell}(f) = k$. In this case, we have $w_{\ell}(f - \theta^{\ell - 1}f) \leq k$. Furthermore, we find that
\begin{equation}\label{theta-and-Ul}
f - \theta^{\ell - 1}f = f \mid U_{\ell} \mid V_{\ell} \equiv (f \mid U_{\ell})^{\ell} \pmod{\ell}.
\end{equation}
It follows that $w_{\ell}(f - \theta^{\ell - 1}f) = w_{\ell}(f \mid U_{\ell})^{\ell} = \ell \: w_{\ell}(f \mid U_{\ell}) \leq k$. Thus, we have $w_{\ell}(f \mid U_{\ell}) \leq \frac{k}{\ell} < 1$, which in turn implies that $w_{\ell}(f|U_\ell) \in \{ - \infty, 0 \}$. Since $f$ is a cusp form, and so is $f \mid U_{\ell}$. We see that $w_{\ell}(f \mid U_{\ell}) = -\infty$. Hence, we conclude that $f \mid U_{\ell} \equiv 0 \pmod{\ell}$. The other direction follows immediately using \eqref{theta-and-Ul}. 
\end{proof}
We now analyze the $\theta$-cycle $\{w_{\ell}(\theta f), \ldots,  w_{\ell}(\theta^{\ell - 1} f)   \}$ in the range $k < \ell < 2k$. 
\begin{proposition} \label{min_j}
Let $k < \ell < 2k$. Then the following holds
\begin{enumerate}
    \item There exists a minimal $j \geq 1$ with $w_{\ell}(\theta^j f) \equiv 0 \pmod{\ell}$, and it equals $j_{1} = \ell - k$. Furthermore, there exists $1 \leq s_{1} \leq \ell - k + 2$ such that $w_{\ell}(\theta^{\ell - k + 1} f) = \ell - k + 3 + ((\ell - k + 2) - s_{1})(\ell - 1)$.
    \item There exists a minimal $j > j_{1}$ with $w_{\ell}(\theta^{j} f) \equiv 0 \pmod{\ell}$, and it has $j = j_{2} \in \{\ell - 1, \ell - 2\}$. Moreover, if $j_2 = \ell - 2$, then $s_1 = \ell - k + 2$, and if $j_2 = \ell - 1$, then we have $s_1 = \ell - k + 1$.
\end{enumerate}
\end{proposition}
\begin{proof}
We begin with the proof of part (1) of Proposition \ref{min_j}. We use induction on $j$ to show that for all $0 \leq j \leq \ell - k - 1$, we have $w_{\ell}(\theta^j f) = k + j(\ell + 1) \not\equiv 0 \pmod{\ell}$. Thus, it follows that $w_{\ell}(\theta^{\ell - j} f) = k + (\ell - k)(\ell + 1) \equiv 0 \pmod{\ell}$. Hence,  $j_{1} = \ell - k$ exists, and further, there exists $s_1 \geq 1$ such that
\begin{align*}
    w_{\ell}(\theta^{\ell - k + 1} f) = k + (\ell - k + 1)(\ell + 1) - s_1 (\ell - 1) = \ell - k + 3 + ((\ell - k + 2) - s_1)(\ell - 1).
\end{align*}
Finally, to show that $s_{1} \leq \ell - k + 2$, we argue by contradiction, using that $f$ is a cusp form with $\ell > k > 3$, and that $\theta^{\ell} f \equiv \theta f \pmod{\ell}$. 
\par
We proceed to the proof of part (2). Suppose that a minimal $j_2$ does not exist. Then for all $1 \leq i \leq k - 1$, we have $w_{\ell}(\theta^{j_1 + i} f) = w_{\ell}(\theta^{\ell - k + i} f) \not\equiv 0 \pmod{\ell}$.
Using part (1), it follows that $w_{\ell}(\theta^{\ell} f) = w_{\ell}(\theta^{\ell - 1} f) + \ell + 1 = k + \ell(\ell + 1)  - s_{1}(\ell - 1)$ and using $w_{\ell}(\theta^{\ell} f) =  w_{\ell}(\theta f) = k + \ell + 1$, we deduce that $s_{1} = \ell + 1$. This leads to a contradiction since $1 \leq s_1 \leq \ell - k + 2 < \ell + 1$. Hence, $j_2$ exists, and there exists some $1 \leq i \leq k - 1$ such that $j_2 = j_1 + i  = \ell - k + i$. Therefore, we have, $w_{\ell}(\theta^{j_2} f) = k + (\ell - k + i) (\ell + 1) - s_{1} (\ell - 1) \equiv 0 \pmod{\ell}$ which implies that $i + s_1 \equiv 0 \pmod{\ell}$. Since $1 \leq i \leq k - 1, 1 \leq s_1 \leq \ell - k + 2$ and $k <  \ell$, we deduce that $i + s_1 = \ell$. Hence, we conclude that $j_2 = j_1 + i \in \{\ell - 2, \ell - 1\}$ and furthermore, we have
\[
s_{1} = \ell - i = \begin{cases}
    \ell - k + 2, \quad \text{if} \; j_2 = \ell - 2 \\
    \ell - k + 1, \quad \text{if} \; j_1 = \ell - 1. \\
\end{cases}
\]
\end{proof}
\begin{lemma}\label{j2_s1_values}
 Suppose $k < \ell < 2k$. Then exactly one of the following holds:
 \begin{enumerate}
     \item If $f \mid U_{\ell} \equiv 0 \pmod{\ell}$, then $(j_2, s_1) = (\ell - 2, \ell - k + 2)$. 
     \item If $f \mid U_{\ell} \not\equiv 0 \pmod{\ell}$, then $(j_2, s_1) = (\ell - 1, \ell - k + 1)$.
 \end{enumerate}
\end{lemma}
\begin{proof}
    For part (1), we let $f \mid U_{\ell} \equiv 0 \pmod{\ell}$. Proposition \ref{Ul_filtration} implies that $w_{\ell}(\theta^{\ell - 1} f) = w_{\ell}(f)$. Suppose that $j_2 \neq \ell - 2$. Using part (2) of Proposition \ref{min_j}, it follows that $j_2 = \ell - 1$ and $s_1 = \ell - k + 1$. Then we have
    \begin{equation}\label{fil_match}
    w_{\ell}(\theta^{\ell - 1} f) = k + (\ell - 1)(\ell + 1) - (\ell - k + 1)(\ell - 1) = k = w_{\ell}(f).
    \end{equation}
On simplifying \eqref{fil_match}, we obtain $k = 0$, leading to a contradiction since $w_{\ell}(f) = k \not\equiv 0 \pmod{\ell}$. \\ For part (2), we let $f \mid U_{\ell} \not\equiv 0 \pmod{\ell}$. Proposition \ref{Ul_filtration} implies that $w_{\ell}(\theta^{\ell - 1} f) \equiv 0 \pmod{\ell}$. Since $f$ is a cusp form, we have $w_{\ell}(\theta^{\ell - 1} f) \geq 1$. Suppose that $j_2 \neq \ell - 1$. Using Proposition \ref{min_j}, it follows that $j_2 = \ell - 2$ and $s_1 = \ell - k + 2$, and hence $w_{\ell}(\theta^{\ell - 2} f) \equiv 0 \pmod{\ell}$. Furthermore,
\[
w_{\ell}(\theta^{\ell - 1} f) = k + (\ell - 1)(\ell + 1) - (\ell - k + 2)(\ell - 1) - s_{2}(\ell - 1) \geq 1
\]
implies that there exists $1 \leq s_2 \leq k - 1$ such that $w_{\ell}(\theta^{\ell - 1} f) = k \ell - (\ell - 1)(s_2 + 1) \equiv s_2 + 1 \pmod{\ell}$. We claim that $s_2 = k - 1$. Since $k < \ell$, we note that $w_{\ell}(\theta^{\ell - 1} f) \equiv s_2 + 1 \not\equiv 0 \pmod{\ell}$. Further, using $w_{\ell}(\theta^{\ell} f) = w_{\ell}(\theta^{\ell - 1} f) + \ell + 1 = w_{\ell}(\theta f)$, we obtain $s_2 = k - 1$ which implies that $w_{\ell}(\theta^{\ell - 1} f) = k$, contradicting $w_{\ell}(\theta^{\ell - 1} f) \neq w_{\ell}(f) = k$. 
\end{proof}

\begin{figure}[H]
\centering
\begin{subfigure}{0.5 \textwidth}
\centering
\begin{tikzpicture}[
    x = 1cm,y = 1cm,
    axis/.style={thick,-{Stealth[length=2mm]}},
    curve/.style={thick,draw=blue!100!black},
    tick/.style={thick},
    every node/.style={font=\small}
  ]

  %--- axes
  \draw[axis] (0,0) -- (12,0) node[below right]{$j$};
  \draw[axis] (0,0) -- (0,8) node[above left]{$w_{\ell}(\Theta^{j}f)$};

  \node[left] at (0,1) {$k$};

  \foreach \x/\lab in {0/{0}, 1/{1}, 3/{\ell-k}, 4.5/{\ell-k+1},
                      7.2/{\ell-2}, 8.8/{\ell-1}, 10.4/{\ell}}{
    \draw[tick] (\x,0.10) -- (\x,-0.10);
    \node[below] at (\x,0) {$\lab$};
  }

  \coordinate (A) at (0,1);   % j=0, value k
  \coordinate (B) at (1,2.5);   % j=1, value k+l +1
  \coordinate (C) at (3,5.4);   % j=l-k, value l(l-k+1) 
  \coordinate (D) at (4.6,1);   % j=l-k+1, value l-k+3    
  \coordinate (E) at (7.2,5.4);   % j=l-2, value l(k-2)   
  \coordinate (F) at (8.8,1);   % j=l-1, value k        
  \coordinate (G) at (10.4,2.5);  % j=l, value k+l+1

  \draw[curve] (A)--(B)--(C)--(D)--(E)--(F)--(G);

   \foreach \P in {A,B,C,D,E,F, G}
    \fill[blue!100!black] (\P) circle (1.6pt);
    
  \node[above right]at (B) {$k+\ell+1$};
  \node[above]      at (C) {$\ell(\ell-k+1)$};
  \node[below=2pt]  at (D) {$\ell-k+3$};
  \node[above]      at (E) {$\ell(k-2)$};
  \node[below right]at (F) {$k$};
  \node[above right]at (G) {$k+\ell+1$};

\end{tikzpicture}
\caption{$f \mid U_{\ell}\equiv 0 \pmod{\ell}
  \iff$ \\ $w_{\ell}(\theta^{\ell - 1}f) = w_{\ell}(f) = k$.}\label{fig1}
\end{subfigure}%
\begin{subfigure}{0.5 \textwidth}
\centering
\begin{tikzpicture}[
    x = 1cm,y = 1cm,
    axis/.style={thick,-{Stealth[length=2mm]}},
    curve/.style={thick,draw=blue!100!black},
    tick/.style={thick},
    every node/.style={font=\small}
  ]

  %--- axes
  \draw[axis] (0,0) -- (12,0) node[below right]{$j$};
  \draw[axis] (0,0) -- (0,8) node[above left]{$w_{\ell}(\Theta^{j}f)$};

  \node[left] at (0,1) {$k$};

  \foreach \x/\lab in {0/{0}, 1/{1}, 3/{\ell-k}, 4.5/{\ell-k+1},
                      7.2/{\ell-2}, 8.8/{\ell-1}, 10.4/{\ell}}{
    \draw[tick] (\x,0.10) -- (\x,-0.10);
    \node[below] at (\x,0) {$\lab$};
  }

  \coordinate (A) at (0,1);   % j=0, value k
  \coordinate (B) at (1,2.4);   % j=1, value k+l +1
  \coordinate (C) at (3,5.4);   % j=l-k, value l(l-k+1) 
  \coordinate (D) at (4.6,2.5);   % j=l-k+1, value l-k+3    
  \coordinate (E) at (7.2,5.7);   % j=l-2, value l(k-2)   
  \coordinate (F) at (8.8,7.5);   % j=l-1, value k        
  \coordinate (G) at (10.4,2.5);  % j=l, value k+l+1

  \draw[curve] (A)--(B)--(C)--(D)--(E)--(F)--(G);

   \foreach \P in {A,B,C,D,E,F, G}
    \fill[blue!100!black] (\P) circle (1.6pt);
    
  \node[above right]at (B) {$k+\ell+1$};
  \node[above]      at (C) {$\ell(\ell-k+1)$};
  \node[below=2pt]  at (D) {$2 \ell-k+2$};
  \node[above]      at (E) {$k (\ell - 1) - 1$};
  \node[below right]at (F) {$k\ell$};
  \node[above right]at (G) {$k+\ell+1$};

\end{tikzpicture}
\caption{$f \mid U_{\ell} \not\equiv 0 \pmod{\ell}
  \iff$ \\ $w_{\ell}(\theta^{\ell - 1}f) \equiv 0 \pmod{\ell}$.} \label{fig2}
\end{subfigure}
\end{figure}

\iffalse
As a corollary, we obtain the following:
\begin{enumerate}
    \item If $f \mid U_{\ell} \equiv 0 \pmod{\ell}$, then
    \[
     (j_1, s_1) = (\ell - k, \ell - k + 2), \qquad
        (j_2, s_2) = (\ell - 2, k - 1), \qquad
 w_{\ell}(\theta^{\ell - k + 1} f) = \ell - k + 3.
 \]
\item If $f \mid U_{\ell} \not\equiv 0 \pmod{\ell}$ then  
\[
(j_1, s_1) = (\ell - k, \ell - k + 1), \qquad (j_2, s_2) = (\ell - 1, k), \qquad
 w_{\ell}(\theta^{\ell - k + 1} f) = 2 \ell - k + 2.
 \]
\end{enumerate}
\fi
We now proceed to show that if $k < \ell < 2k$, then we have $\ell \in \{2k - 3, 2k - 1\}$.
\iffalse
\begin{theorem}
Let $\ell > 3$ be prime. Let $f$ be as defined above with a Type II congruence modulo $\ell$, and let $k < \ell < 2k$. Then exactly one of the following holds:
\begin{enumerate}
    \item If $f \mid U_{\ell} \equiv 0 \pmod{\ell}$, then we have $\ell = 2k - 3$.
    \item If $f \mid U_{\ell} \not\equiv 0 \pmod{\ell}$, then we have  $\ell = 2k - 1$.
\end{enumerate}
\end{theorem}
\fi
 Let $\ell > 3$ be exceptional of Type II with $k < \ell < 2k$. We observe that $\ell - k + 1 \leq \frac{\ell + 1}{2} \leq \ell - 2$. Using \eqref{type2_form}, we have 
 \[
 w_{\ell}(\theta (f \otimes 1_N)) =  w_{\ell}(\theta f) = w_{\ell}(\theta^{\frac{\ell + 1}{2}} f) = w_{\ell}(\theta^{\frac{\ell + 1}{2}} (f \otimes 1_N)).
 \]
 Furthermore, Lemma \ref{j2_s1_values}, figures \eqref{fig1} and \eqref{fig2} implies that
 \begin{equation}\label{fil_comparison}
k + \ell + 1 = 
 \begin{cases}
     k + \frac{\ell + 1}{2}(\ell + 1) - (\ell - k + 2)(\ell - 1),  \:\: f \mid U_{\ell} \equiv 0 \pmod{\ell}, \\
     k + \frac{\ell + 1}{2}(\ell + 1) - (\ell - k + 1)(\ell - 1),  \:\: f \mid U_{\ell} \not\equiv 0 \pmod{\ell}.
 \end{cases}
 \end{equation}
 On simplifying \eqref{fil_comparison}, we find that when $f \mid U_{\ell} \equiv 0 \pmod{\ell}$, we have $\ell = 2k - 3$, while when $f \mid U_{\ell} \not\equiv 0 \pmod{\ell}$, we have $\ell = 2k - 1$.

\bibliographystyle{alpha}  
\bibliography{bib}
\end{document}